\documentclass[11pt,reqno]{amsart}
\usepackage{amssymb,amsthm,amsmath,comment}

\begin{document}

\newtheorem{thm}{Theorem}[section]
\newtheorem{cor}[thm]{Corollary}
\newtheorem{lem}[thm]{Lemma}
\newtheorem{counterex}[thm]{Counterexample}
\newtheorem{defandth}[thm]{Definition and Theorem}
\newtheorem*{conj}{Conjecture}
\theoremstyle{definition}
\newtheorem{defn}[thm]{Definition}
\theoremstyle{remark}
\newtheorem{rem}[thm]{Remark}
\newcommand{\xvec}{x_1,\ldots,x_k}
\newcommand{\yvec}{y_1,\ldots,y_\ell}
\newcommand{\im}{{\rm im}}
\newcommand{\I}{\mathcal{I}}
\newcommand{\ds}{doubly symmetric}

\title{Doubly Symmetric Functions}

\author[Allan Berele]{Allan Berele$^\dagger$}
\author[Bridget Eileen Tenner]{Bridget Eileen Tenner}
\thanks{$^\dagger$ Work supported by the National Security Agency, under Grant H98230-08-1-0026.} 
\email{aberele@condor.depaul.edu}\email{bridget@math.depaul.edu}
\address{Department of Mathematical Sciences, DePaul University, Chicago, Illinois 60614}


\begin{abstract}
In this paper we introduce doubly symmetric functions, arising from the equivalence of particular linear combinations of Schur functions and hook Schur functions.  We study algebraic and combinatorial aspects of doubly symmetric functions, in particular as they form a subalgebra of the algebra of symmetric functions.  This subalgebra is generated by the odd power sum symmetric functions.  One consequence is that a Schur function itself is doubly symmetric if and only if it is the Schur function of a staircase shape.
\end{abstract}

\maketitle

It is well-known that the Schur functions $S_\lambda$ span the
ring of symmetric functions.  These were generalized to the hook
Schur functions $HS_\lambda$ in \cite{BR}, and both are indexed by
partitions.  In this paper we define a new class of symmetric functions, defined by the coincidence of linear combinations of Schur and hook Schur functions.  More precisely, a symmetric function $\sum m_\lambda S_\lambda$ is {\it doubly symmetric\/} if for all (possibly empty) sets of variables $x_1,\ldots,x_k$ and $y_1,\ldots,y_\ell$,
$$\sum m_\lambda S_\lambda(\xvec,\yvec) = \sum m_\lambda HS_\lambda(x_1,\ldots,x_k;y_1,\ldots,y_\ell).$$
  One example of a doubly symmetric function is the Schur function $S_{(1)}$, since
$$S_{(1)}(x_1,\ldots,x_n)=x_1+\cdots+x_n=HS_{(1)}(x_1,\ldots,x_i;x_{i+1},\dots,x_n).$$
A more interesting example is $k_n=\sum_{a+b=n} S_{(a,1^b)}$.  That
$$\sum_{a+b=n} S_{(a,1^b)}(\xvec,\yvec) =\sum_{a+b=n} HS_{(a,1^b)}(\xvec;\yvec)$$
can be proven bijectively as follows.  We denote by $H(1,1)$ the set of partitions of the form
$(a,1^b)$.  That is, elements of $H(1,1)$ are the partitions consisting of at most one row of length greater than 1 and at most one column of length greater than 1.  Given an ordinary semistandard tableau of shape $\lambda \in H(1,1)$ filled by $k+\ell$ letters, we produce a $(k,\ell)$-semistandard tableau, defined in Section~\ref{sec:hook schur}, by moving the letters greater than $k$ from the first column to the first row, and {\it vice versa}.

It turns out that the product of doubly symmetric functions is
again doubly symmetric, and so they form a subalgebra
$\mathcal{D}$ of the symmetric functions~$\mathcal{S}$.  There is
an onto homomorphism $D:\mathcal{S}\rightarrow\mathcal{D}$ defined
on Schur functions via
$$D(S_\lambda)=\sum_{\mu\subseteq \lambda} S_\mu S_{(\lambda/\mu)'}.$$
We denote $D(S_\lambda)$ by $DS_\lambda$ and call it a doubly Schur function.  The functions
$\{DS_\lambda\}$ have many properties in common with ordinary
Schur functions and many Schur function identities have ``double'' analogues.  These functions span $\mathcal{D}$ but are not linearly
independent.  The kernel of~$D$ is the ideal~$\I$ generated by all
$S_\lambda-S_{\lambda'}$, and so
$$\mathcal{D}\simeq \mathcal{S}/\I.$$
Moreover, a stronger property is true:  $\mathcal{D}$ is the orthogonal compliment of~$\I$ with respect to the natural inner product in the ring of symmetric functions and so
$\mathcal{S}=\mathcal{D}\oplus \I$.

Recall the power symmetric functions, defined by $p_n=\sum_i x_i^n$ and
$p_\lambda=p_{\lambda_1}p_{\lambda_2}\cdots p_{\lambda_m}$ for a partition $\lambda$ with parts $\lambda_1 \ge \lambda_2 \ge \cdots \ge \lambda_m > 0$.  These, like the Schur
functions, form a basis of~$\mathcal{S}$.  In fact, we show that each $p_\lambda$ is an
eigenvector for $D$.  If some part $\lambda_i$ is even, then
$D(p_\lambda)=0$; otherwise, $D(p_\lambda)=2^mp_\lambda$, where $m$ is
the length of~$\lambda$.  This implies that $\{p_\lambda \mid \lambda
\mbox{ has all odd parts}\}$ is a basis for $\mathcal{D}$ and, in
particular, the dimension of the space of doubly symmetric
functions of degree~$n$ is the number of partitions of~$n$ into
odd parts; that is, the coefficient of $x^n$ in the generating function
$$\prod_{k\ge0} (1-x^{2k+1})^{-1}.$$
It is well-known that the number of partitions of $n$ into odd parts is equal to the number of partitions of $n$ into distinct parts.  Thus, we can also say that the dimension of the space of doubly symmetric functions of degree~$n$ is the coefficient of $x^n$ in the generating function
$$\prod_{k \ge 1} (1 + x^k).$$

In the first section of this paper, we review the properties and results of hook Schur functions which are relevant to this work.  In Section~\ref{section:doubly}, we formally introduce doubly symmetric functions, deriving both algebraic and combinatorial results for these functions in general and for doubly Schur functions in particular.  The final section of the paper contains the main results, including a description of the relationship between the algebra of symmetric functions and the subalgebra of doubly symmetric functions (Theorem~\ref{thm:subalgebra}).  It follows from this result that the subalgebra of doubly symmetric functions is generated by the odd power sum symmetric functions, and moreover that a Schur function $S_{\lambda}$ itself is doubly symmetric if and only if $\lambda$ is a staircase shape.

\section{Hook Schur Functions}\label{sec:hook schur}

We assume that the reader is familiar with the theory of Schur
functions, and we review the theory of
hook Schur functions.  Except as indicated, this material can be
found in \cite{BR}.  More information about Schur functions can be found, for example, in \cite[Chapter I]{M}.

\begin{defn}
Given  an ordered alphabet with $k$ ordinary symbols
and $\ell$ primed symbols -- in this paper we will use the alphabet
$\{1<2<\cdots< k< 1'<\cdots< \ell'\}$ -- and a partition
$\lambda$, we defined a \emph{$(k,\ell)$-semistandard tableau of shape
$\lambda$} to be a tableau of shape $\lambda$ subject to:
\begin{enumerate}
\item entries are non-decreasing in each row and column, and
\item unprimed symbols may repeat in rows, but not columns, while
prime symbols may repeat in column, but not rows.
\end{enumerate}
\end{defn}

Given such a tableau we write $x^T$ for the monomial in
$x_1,\ldots,y_{\ell}$ whose degree in each $x_j$ equals the
number of~$j$ in~$T$, and whose degree in each $y_j$ is the number of $j'$.  See Figure~\ref{fig:(2,3)-sst} for a (2,3)-semistandard tableau $T$ of shape $(5,4,2,2,1,1)$, with
$x^T=x_1^3x_2^3 y_1^4 y_2^4 y_3$.

\begin{figure}[htbp]
$$\begin{array}{ccccc}
1&1&1&2&3'\\ 2&2&1'&2'&\\ 1'&2'&&\\ 1'&2'&&\\ 1'&&&&\\2'&&&&
\end{array}$$
\caption{A (2,3)-semistandard tableau of shape $(5,4,2,2,1,1)$.}\label{fig:(2,3)-sst}
\end{figure}

\begin{defn}
For a partition $\lambda$, the \emph{hook Schur function}
$$HS_\lambda(\xvec;\yvec)$$
is equal to $\sum \{x^T \mid T \mbox{ is a $(k,\ell)$-semistandard tableau of shape }
\lambda\}.$  We will sometimes abbreviate this as
$HS_\lambda(X;Y)$, where the $X$ and $Y$ stand for sets of
variables.  Skew hook Schur functions $HS_{\lambda/\mu}$ are
defined analogously using skew tableaux.
\end{defn}

This formula is equivalent to the definition (compare to \cite[\S I.5 (5.10)]{M}):
\begin{equation}\label{eq:1}
HS_\lambda(X;Y)=\sum_{\mu\subseteq \lambda} S_\mu(X)S_{(\lambda/\mu)'}(Y).
\end{equation}
One consequence, which follows from the corresponding formula for Schur functions, is that for fixed $\lambda$
\begin{equation}\label{eq:1a}
\sum_{\mu\subseteq\lambda} HS_\mu(X;Y)HS_{\lambda/\mu}(Z;U)=HS_\lambda(X,Z;Y,U).
\end{equation}
A less obvious formula we will need is
\begin{equation}\label{eq:2}
HS_\lambda(X;Y)=HS_{\lambda'}(Y;X).
\end{equation}
The symmetry of the $HS_\lambda$ in the next theorem was first proven by
Berele and Regev, see~\cite{BR}, and a combinatorial proof due to Remmel can be found in~\cite{R}. The rest of Theorem~\ref{thm:symmetry} is due to Stembridge in~\cite{S}.

\begin{thm}\label{thm:symmetry}
For all $\lambda$, the hook Schur function $HS_\lambda(x_1,\ldots,x_k;\yvec)$ is
symmetric in $x_1,\ldots,x_k$ and in $\yvec$,
and has the additional property that for each $a\le k$, $b\le
\ell$, $$HS_\lambda(\xvec;\yvec)\vert_{x_a=-y_b=t}$$ is
independent of~$t$.  Moreover, the hook Schur functions span the
ring of functions with these properties.
\end{thm}

There are analogues of the Robinson-Schensted correspondence for hook Schur functions, see~\cite{BRl}.  They imply that 
\begin{equation}
\sum_{\lambda\in{\rm Par}(n)} f^\lambda HS_\lambda(\xvec;\yvec)=(x_1+\cdots+y_\ell)^n,\label{eq:2a}
\end{equation}
 where $f^\lambda$ equals the number of standard Young tableaux of shape~$\lambda$; and an analogue of the Cauchy identity:
\begin{multline}\label{eq:3}
\sum_\lambda HS_\lambda(x_1,\ldots,x_a;y_1,\ldots,y_b)HS_\lambda(z_1,\ldots,z_c;
u_1,\ldots,u_d)\\
=\prod(1+x_\alpha u_\delta)\prod(1+y_\beta
z_\gamma)\prod(1-x_\alpha z_\gamma)^{-1}\prod(1-y_\beta
u_\delta)^{-1},
\end{multline}
where the products are indexed by $\alpha \in \{1,\ldots,a\}$,  $\beta \in \{1,\ldots,b\}$, and so on.

\begin{defn}
$H(k,\ell;n)$ is the set of  partitions of~$n$ with at
most $k$ parts greater than~$\ell$, and $H(k,\ell)$ is the
union $\cup_n H(k,\ell;n)$.
\end{defn}

In other words, no shape in $H(k, \ell)$ contains a rectangle of height $k+1$ and width $\ell+1$.

\begin{thm}
$HS_\lambda(x_1,\ldots,x_k;\yvec)\ne0$ if and only if $\lambda\in H(k,\ell)$.
\end{thm}

\begin{defn}
Let $\mathcal{H}(k,\ell)$ be the ring generated by
the hook Schur functions
$$HS_\lambda(\xvec;\yvec).$$
There are maps from $\mathcal{H}(k,\ell)$ to each of $\mathcal{H}(k-1,\ell)$ and $\mathcal{H}(k,\ell-1)$ obtained by setting $x_k$ or $y_\ell$ equal to zero, and these maps take $H_\lambda(\xvec;\yvec)$ to $HS_\lambda(x_1,\ldots,x_{k-1};\yvec)$ and $HS_\lambda(\xvec;y_1,\ldots,y_{\ell-1})$, respectively.  This allows us to define $\mathcal{H}$ as the inverse limit of these rings, and $HS_\lambda$ as the element of $\mathcal{H}$ corresponding to the $HS_\lambda(X;Y)$.
\end{defn}

\begin{thm}\label{thm:1}
The linear map
$\digamma:\mathcal{S}\rightarrow\mathcal{H}$ given by
$\digamma(S_\lambda)=HS_\lambda$ is an algebra
homomorphism.
\end{thm}

(The archaic Greek $\digamma$ is a wau, which means hook. Thanks are due to the first author's wife Miriam for suggesting this notation.)

We now compute the image of the power symmetric function $p_n$
under $\digamma$.  We note that these functions are discussed by Kantor and Trishin in~\cite{KT}.

\begin{lem}
If $n\ge1$ then $p_n=\sum_{a=1}^n
(-1)^{n-a}S_{(a,1^{n-a})}$.\label{lem:MN}
\end{lem}

\begin{proof}
By the Murnaghan-Nakayama rule (see \cite[\S I.7]{M}), $p_n=\sum_\lambda (-1)^{ht(\lambda)}S_\lambda$, where the sum is over all border strips having size~$n$.  The only such shapes are of the form $(a,1^{n-a})$, and the height of this shape is $n-a$.
\end{proof}

\begin{thm}\label{thm:p}
$\digamma(p_n)(x_1,\ldots,x_k;\yvec)$ equals
$$x_1^n+\cdots+x_k^n+(-1)^{n+1}(y_1^n+\cdots+y_\ell^n).$$
\end{thm}

\begin{proof}
By the previous lemma we know that
\begin{equation}\label{eq:p_n computation}
\digamma(p_n)(x_1,\ldots,x_k;\yvec) =\sum_{a=1}^n (-1)^{n-a} HS_{(a,1^{n-a})}(\xvec;\yvec).
\end{equation}
Any $(k,\ell)$-semistandard tableau of shape $\lambda \in H(1,1)$ which has two different variables in its filling corresponds bijectively to a filling of the same content of some $\mu \in H(1,1)$, where $|ht(\lambda) - ht(\mu)| = 1$.  Hence these two terms would cancel in the sum.

Suppose that there is such a shape and filling containing both the unprimed numbers $j>i$.  Since $j$ is not the smallest letter occurring in the tableau, it cannot occupy the unique box which is in both the top row and the left column of the shape.  If $j$ appears outside the top row (that is, if it occurs in the left column), then move the box for that appearance to the top row.  If $j$ does not appear outside the top row, then move a~$j$ box from the top row to the left column.  These operations are inverses, change the height of the shape by~1, and neither changes the content of the filling.  Thus, any such terms cancel in the sum on the right side of equation~\eqref{eq:p_n computation}.  Tableaux using two different primed numbers can be handled analogously (swapping the roles of the left column and top row).  Finally, we consider shapes filled entirely by some~$i$ and $j'$.  If the $j'$ occurs (once) in the top row, then move that box to the left column.  Otherwise, do the reverse procedure.  As in the previous cases, this again causes the terms to cancel.

Thus we are left with only the fillings involving only a single letter.  If the letter is unprimed, then the shape must be~$(n)$, having height~0.  If the letter is primed then the shape is $(1^n)$, having height~$n-1$.  This completes the proof.
\end{proof}

\section{Doubly Symmetric Functions}\label{section:doubly}

Recall the following definition, mentioned in the introduction.

\begin{defn}
The function $\sum m_\lambda S_\lambda$ is
\emph{doubly symmetric} if, for all finite sets of variables $X$
and $Y$,
\begin{equation}\label{eq:5}
\sum m_\lambda S_\lambda(X,Y)= \sum m_\lambda HS_\lambda(X;Y).
\end{equation}
The set of all doubly symmetric functions is denoted $\mathcal{D}$.
\end{defn}

The following lemma describes elementary properties of these functions.

\begin{lem}
\begin{enumerate}
\item If $\sum m_\lambda S_\lambda$ is doubly symmetric, then 
so is the sum $\sum\{ m_\lambda S_\lambda \mid \lambda\in
\mbox{Par}(n)\}$ for each $n$.
\item If $\sum m_\lambda S_\lambda$ is doubly symmetric then
$m_\lambda=m_{\lambda'}$ for each~$\lambda$.
\item The product of doubly symmetric functions is doubly symmetric.
\end{enumerate}
\end{lem}

\begin{proof}
\begin{enumerate}
\item This is obvious.
\item This follows from equation (\ref{eq:5}), taking $X=\emptyset$.
\item This follows from Theorem~\ref{thm:1}.
\end{enumerate}
\end{proof}

\begin{defn}
Let $\I$ be the ideal of $\mathcal{S}$ generated by all differences $S_\lambda-S_{\lambda'}$.
\end{defn}

\begin{thm}
The subalgebra $\mathcal{D}$ and the ideal $\I$ are orthogonal compliments with respect to the inner product in
$\mathcal{S}$.\label{thm:9}
\end{thm}

\begin{proof}
Let $C_{\mu,\nu}^\lambda$ be the
Littlewood-Richardson coefficients, so
$$S_\mu\cdot S_\nu =\sum_\lambda C_{\mu,\nu}^\lambda S_\lambda\mbox{ and
}S_{\lambda/\mu}=\sum_\nu C_{\mu,\nu}^\lambda S_\nu.$$
Equivalently, $C^\lambda_{\mu,\nu}=\langle S_\lambda,S_\mu S_\nu\rangle$. Then in
equation~\eqref{eq:5} we first apply \cite[\S I.5 (5.10)]{M} to the
left hand side to get
\begin{eqnarray*}
\sum_\lambda m_\lambda S_\lambda(X,Y)&=&
\sum_\lambda m_\lambda \sum S_\mu(X)S_{\lambda/\mu}(Y)\\
&=& \sum_\lambda m_\lambda\sum_{\mu,\nu}C_{\mu,\nu}^\lambda
S_\mu(X)S_\nu(Y).
\end{eqnarray*}
Likewise, applying equation~\eqref{eq:1} to the left
hand side of equation~\eqref{eq:5} gives
$$\sum_\lambda m_\lambda
\sum_{\mu,\nu}C_{\mu',\nu}^{\lambda'}S_\mu(X) S_\nu(Y).$$
We now equate the coefficients of $S_\mu(X)S_\nu(Y)$ to get
$$\sum_\lambda m_\lambda C_{\mu,\nu}^\lambda =\sum_\lambda
m_\lambda C_{\mu',\nu}^{\lambda'}.$$
Taking into account $C_{\mu',\nu}^{\lambda'}=C_{\mu,\nu'}^\lambda$ we get
$$0=\sum_\lambda m_\lambda(C_{\mu,\nu}^\lambda-C_{\mu,\nu'}^\lambda)=\left\langle
\sum_\lambda m_\lambda S_\lambda, S_\mu(S_\nu-S_{\nu'})\right\rangle.$$
Hence, the sum $\sum m_\lambda S_\lambda$ is doubly symmetric if and only
if it is orthogonal to the ideal generated by all $S_\nu-S_{\nu'}$.
\end{proof}

\begin{cor}
The ring of doubly symmetric functions is isomorphic
to the ring of symmetric functions, modulo the ideal $\I$.
\end{cor}

\begin{defn}
Let $D:\mathcal{S}\rightarrow\mathcal{S}$ be the
linear map defined by $D(S_\lambda)=\sum_\mu S_\mu
S_{(\lambda/\mu)'},$ and extended to all of $\mathcal{S}$ by
linearity. We write $D(S_\lambda)$ as $DS_\lambda$ and call it a
\emph{doubly Schur function}.
\end{defn}

Lemma~\ref{lem:properties of doubly schur} lists some basic
properties of doubly Schur functions, including the fact that they
are doubly symmetric.

\begin{defn}
As in the opening discussion of this article, set $k_n = \sum_{a+b=n}S_{(a,1^b)}$.
\end{defn}

\begin{lem}\label{lem:properties of doubly schur}
\begin{enumerate}
\item $DS_{(n)}=2k_n$.
\item $DS_\lambda(x_1,\ldots,x_n)=HS_\lambda(x_1,\ldots,x_n;x_1,\ldots,x_n)$.
Hence, by linearity, $Df(X)=\digamma(f)(X;X)$ for all
$f\in\mathcal{S}$.
\item $DS_\lambda$ is doubly symmetric.  Hence, by linearity,
$D(\mathcal{S})\subseteq\mathcal{D}$.
\item $D$ is an algebra homomorphism.
\item If $*$ is the internal product of symmetric functions and $f$ is a degree~$n$ symmetric function, then $D(f)=f*2k_n$.
\end{enumerate}
\end{lem}

\begin{proof}
\begin{enumerate}
\item We compute
$$DS_{(n)} = \sum_{i=0}^n
S_{(i)}S_{(n-i)'}=\sum_{i=0}^n S_{(i)} S_{(1^{n-i})}.$$
If $1\le i\le n-1$, then
$S_{(i)}S_{(1^{n-i})}=S_{(i,1^{n-i})}+S_{(i+1,1^{n-i-1})}$ and the result
follows.
\item This follows from comparing $DS_\lambda(X)=\sum_\mu
S_\mu(X)S_{(\lambda/\mu)'}(X)$ with $HS_\lambda (X;Y)=\sum_\mu
S_\mu(X) S_{(\lambda/\mu)'}(Y),$ in which we set $X=Y$.
\item Let $\digamma$ be as in Theorem~\ref{thm:1}.  We consider
$\digamma(DS_\lambda)(X;Y)$ and show that it equals $DS_\lambda(X,Y)$.
$$\begin{array}{lll}
\digamma(DS_\lambda)(X;Y)&= \digamma(\sum_\mu
S_\mu(X,Y) S_{(\lambda/\mu)'}(X,Y))&\mbox{(by definition of }\mathcal{D})\\
&= \sum_\mu HS_\lambda(X;Y) HS_{(\lambda/\mu)'}(X;Y)&\mbox{(by Theorem \ref{thm:1})}\\
&= \sum_\mu HS_\lambda(X;Y) HS_{(\lambda/\mu)}(Y;X)&\mbox{(by equation \eqref{eq:2})}\\
&=HS_\lambda(X,Y;X,Y)&\mbox{(by equation \eqref{eq:1a})}\\
&=DS_\lambda(X,Y)&\mbox{(by part (2) above)}
\end{array}$$
\item This follows from the observation that each of the maps
$$S_\lambda(X)\mapsto HS_\lambda(X;Y) \text{ and } HS_\lambda(X;Y)\mapsto HS_\lambda(X;X)$$
is a homomorphism.
\item This follows from the following computation, using \cite[\S I.7, exercise 23(b)]{M}.
\begin{align*}
DS_\lambda&=\sum_{\mu\subseteq \lambda} S_\mu S_{(\lambda/\mu)'}\\ &=\sum_{i=0}^n \sum_{\substack{\mu\subseteq\lambda\\ |\mu|=i}}S_\mu S_{(\lambda/\mu)'}\\ &=\sum_{i=0}^n \sum_{\substack{\mu\subseteq\lambda\\ |\mu|=i}}(S_\mu*S_{(i)})(S_{\lambda/\mu}*S_{(1^{n-i})})\\ &=\sum_{i=0}^n S_\lambda*(S_{(i)}S_{(1^{n-i})})\\ &=S_\lambda*D(S_{(n)})
\end{align*}
\end{enumerate}
\end{proof}

\begin{rem}\label{rem:1}
In light of the last part of Lemma~\ref{lem:properties of doubly schur}, the doubling map $D$ also has an interpretation in terms of $S_n$-characters.
The Frobenius map is an isomorphism between characters
of the symmetric group $S_n$ and symmetric functions of degree~$n$  such that $Fr(\chi^\lambda)=S_\lambda$.  If $f$ is
 a symmetric function of degree~$n$, then $Fr^{-1}(\mathcal{D}(f))$ is the inner tensor product
$Fr^{-1}(f)\otimes \sum_{a+b=n} 2\chi^{(a,1^b)}.$ 
\end{rem}

It is now straightforward to prove ``double'' analogues of Schur function identities.

\begin{thm}\label{thm:double analogues}
\begin{enumerate}
\item $DS_\lambda =\det(2k_{\lambda_i-i+j})_{i,j}$.
\item $\sum_{\lambda\vdash n} f^\lambda DS_\lambda(\xvec)=4(x_1+\cdots+x_k)^n$
\item If $X=\{x_1,\ldots,x_a\}$ and $Y=\{y_1,\ldots,y_b\}$ then
$$\sum_\lambda DS_\lambda(X) DS_\lambda(Y) = \prod_{i=1}^a
\prod_{j=1}^b \frac{(1+x_iy_j)^2}{(1-x_iy_j)^2}$$ and
\begin{equation}\label{eqn:result of Cauchy}
\sum_\lambda DS_\lambda(X) S_\lambda(Y) = \prod_{i=1}^a
\prod_{j=1}^b \frac{(1+x_iy_j)}{(1-x_iy_j)}.
\end{equation}
\end{enumerate}
\end{thm}

\begin{proof}
\begin{enumerate}
\item This follows from the application of $D$ to the determinantal
identity $S_\lambda = \det (h_{\lambda_i-i+j})$, see \cite[\S I.3 (3.5)]{M}.
\item This follows from equation \eqref{eq:2a}.
\item We apply the hook Schur version of
the Cauchy identity, equation~{\eqref{eq:3}}, to $$\sum_\lambda
DS_\lambda(X) DS_\lambda(Y)=\sum_\lambda HS_\lambda(X;X)
HS_\lambda (Y;Y)$$ and $$\sum_\lambda DS_\lambda(X)
S_\lambda(Y)=\sum_\lambda HS_\lambda(X;X)
HS_\lambda(Y;\emptyset),$$ respectively.
\end{enumerate}
\end{proof}

We may use equation~\eqref{eqn:result of Cauchy} of Theorem~\ref{thm:double analogues} to get an inner product or integral representation
of~$D$.  By Weyl's formula the inner product in $\mathcal{S}$ can
be computed via
$$\langle f,g\rangle = \oint_T f(X)g\left(X^{-1}\right)\, d\nu,$$
where $X=\{x_1,\ldots,x_n\}$, $X^{-1}=\{x_1^{-1},\ldots,x_n^{-1}\}$, $T$ is the torus $|x_i|=1$,
$i=1,\ldots,n$, and
$$d\nu=(2\pi \sqrt{-1})^{-n} (n!)^{-1}\prod_{1\le i\ne j\le
n}\left(1-\frac{x_i}{x_j}\right)\frac{dx_1}{x_1}\wedge\cdots\wedge\frac{dx_n}{x_n}.$$

\begin{thm} \label{thm:14}
If $f(x_1,\ldots,x_n)$ is symmetric, then
\begin{align*}Df(y_1,\dots,y_m)&=\left\langle f,
\prod(1+x_iy_j)(1-x_iy_j)^{-1}\right\rangle \\
&=\oint_T f(x_1,\ldots,x_n)\prod\frac{x_i+y_j}{x_i-y_j} \,d\nu.
\end{align*}
\end{thm}

\begin{proof}
Let $f=\sum m_\lambda S_\lambda$, so
$m_\lambda=\langle f,S_\lambda\rangle$ and $Df=\sum_\lambda
m_\lambda DS_\lambda$.  Hence,
\begin{align*}
Df&=\sum \left\langle
f(X), S_\lambda(X)\right\rangle DS_\lambda(Y) \\
&= \left\langle f(X),\sum
S_\lambda(X)DS_\lambda(Y)\right\rangle\\
&=\left\langle f(X), \prod
(1+x_iy_j)(1-x_iy_j)^{-1}\right\rangle.
\end{align*}
To pass from this to the integral version using Weyl's formula, note that
$$(1+x_i^{-1}y_j)(1-x_i^{-1}y_j)^{-1} = (x_i+y_j)/(x_i-y_j).$$
\end{proof}

Theorem \ref{thm:14} is an analogue for the following formula for $\digamma$, which we have not seen previously.

\begin{thm}
If $f(x_1,\ldots,x_n)$ is symmetric then
\begin{align*}\digamma f(y_1,\dots,y_k;z_1,\ldots,z_\ell)&=\left\langle f,
\prod(1+x_iz_j)(1-x_iy_j)^{-1}\right\rangle \\
&=\oint_T f(x_1,\ldots,x_n)\prod\frac{x_i+z_j}{x_i-y_j}\, d\nu.
\end{align*}
\end{thm}

\section{The Main Theorem}\label{section:main theorem}

\begin{lem}\label{lem:3.1}
$D(p_n)=
\begin{cases}
0 & \text{ if $n$ is even, and}\\
2p_n & \text{ if $n$ is odd.}
\end{cases}$
\end{lem}

\begin{proof}
By Lemma \ref{lem:properties of doubly schur}(2) and Theorem~\ref{thm:p},
\begin{align*}
D(p_n)(x_1,\ldots,x_k)=&\digamma(p_n)(x_1,\ldots,x_k;x_1,\ldots,x_k)\\
=& \sum x_i^n+(-1)^{n+1}\sum x_i^n,
\end{align*}
which equals 0 if $n$ is even, and $2p_n$ if $n$ is odd.
\end{proof}

\begin{defn}
If every part of the partition $\lambda$ is odd, then $\lambda$ is an \emph{odd} partition.
\end{defn}

\begin{cor}\label{cor:1}
Let $\lambda=(\lambda_1,\ldots,\lambda_m)$ be a
partition of length~$m$, and let
$p_\lambda=p_{\lambda_1}\cdots p_{\lambda_m}$ be the corresponding
power symmetric function.  Then $p_\lambda$ is an eigenvector
for~$D$, with
$$D(p_\lambda)=
\begin{cases}
0 & \text{ if some $\lambda_i$ is even, and}\\
2^mp_\lambda & \text{ if $\lambda$ is odd.}
\end{cases}$$
\end{cor}

\begin{cor}\label{cor:2}
The image of $D$ is spanned by the set 
$\{p_\lambda \mid \lambda \text{ is odd}\}$, and the kernel is spanned by the set
$\{p_\lambda \mid \lambda \text{ is not odd}\}$.  In particular,
$\mathcal{S}= \im(D)\oplus \ker(D)$.
\end{cor}

\begin{lem}\label{lem:p_n ideal}
If $n$ is even, then $p_n$ is contained in the ideal $\I$.
\end{lem}

\begin{proof}
The conjugate of the shape $(a,1^b)$ is $(b+1,1^{a-1})$.  Thus, modulo the ideal $\I$, 
we have the congruence $S_{(a,1^b)}\equiv S_{(b+1,1^{a-1})}$.  Letting $n=2r$ and
applying this congruence to Lemma~\ref{lem:MN} yields
\begin{align*}
p_n&=\sum_{a=1}^r\left[(-1)^{n-a}S_{(a,1^{n-a})}+(-1)^{a-1}S_{(n-a+1,1^{a-1})}\right]\\
&\equiv\sum_{a=1}^r\left[(-1)^{n-a}+(-1)^{a-1}\right]S_{(a,1^{n-a})}\\
&=0.
\end{align*}
\end{proof}

\begin{thm}\label{thm:subalgebra}
The map $D:\mathcal{S}\rightarrow\mathcal{D}$ is surjective, with kernel~$\I$.
\end{thm}

\begin{proof} 
By Theorem~\ref{thm:9},
$\mathcal{S}=\mathcal{D}\oplus \I$.  By Corollary~\ref{cor:2},
$\mathcal{S}=\im(D)+\ker(D)$.  The kernel of~$D$ is generated by
$p_{2r}$ and these elements are all contained in~$\I$ by Lemma~\ref{lem:p_n ideal},
so $\ker(D)\subseteq \I$.  By Lemma~\ref{lem:properties of doubly schur}(3),
we have $\im(D)\subseteq\mathcal{D}$.  Hence, $\ker(D)=\I$ and
$\im(D)=\mathcal{D}$.
\end{proof}

\begin{cor}\label{cor:basis of D}
$\mathcal{D}$ is generated by $\{p_n \mid n\mbox{
is odd}\}$ and has basis $\{p_\lambda \mid \lambda\mbox{ is odd}\}$. Also,
because the map $D$ is surjective, $\mathcal{D}$ is spanned by $\{DS_\lambda\}$
and generated by $\{k_n\}$.
\end{cor}

We may also express the $k_n$ in terms of the $p_\lambda$.

\begin{cor}
For each $n\ge 1$, $k_n=\sum 2^{ht(\lambda)-1}z_\lambda^{-1}p_\lambda$, summed over odd partitions of~$n$.
\end{cor}

\begin{proof}
Let $k_n=\sum m_\mu p_\mu$.  We compute $k_n*p_\lambda$ in two different ways.  On the one hand, by \cite[\S I.7 (7.12)]{M}, each $p_\mu*p_\lambda=\delta_{\mu,\lambda}z_\lambda p_\lambda$, and so $k_n*p_\lambda=z_\lambda m_\lambda p_\lambda$.  On the other hand, by part~(5) of Lemma~\ref{lem:properties of doubly schur}, $p_\lambda * k_n=\frac12 D(p_\lambda)$.  By Lemma~\ref{lem:3.1} this equals~0 if $\lambda$ is not odd and  $\frac12 2^{ht(\lambda)}p_\lambda$ if $\lambda$ is odd.   Hence, in this latter case $z_\lambda m_\lambda = 2^{ht(\lambda)-1}$ and the corollary follows.
\end{proof}

Determining the basis of $\mathcal{D}$ allows us to calculate the dimension of this space of polynomials, for any fixed degree, as described in the following corollary.  The sequence of these dimensions is entry A000009 in \cite{oeis}.

\begin{defn}
A partition $\lambda$ is a \emph{staircase} shape if $\lambda = (m,m-1, \ldots, 2, 1)$ for some $m$.
\end{defn}

The following corollary follows from Corollary~\ref{cor:basis of D}, together with Proposition 7.17.7 and Exercise 7.54 of \cite{ec2}.

\begin{cor}\label{cor:staircases}
The Schur function $S_{\lambda}$ is doubly symmetric if and only if $\lambda$ is a staircase shape.
\end{cor}

\begin{cor}
The dimension of the space of doubly symmetric functions of degree $n$ is the coefficient of $x^n$ in
$$\prod_{k\ge0} (1-x^{2k+1})^{-1} = \prod_{k \ge 1} (1+x^k).$$
\end{cor}

Another consequence of Corollary~\ref{cor:basis of D} is the following.

\begin{thm}\label{thm:characterization}
Let $f=\sum m_\lambda S_\lambda$ be such that for each
$0\le i\le n$, the function $f(x_1,\ldots,x_n)$ can be written as a linear
combination
$$\sum c_\lambda HS_\lambda(x_1,\ldots,x_i;x_{i+1},\ldots,x_n),$$
where the $c_\lambda$ potentially depend on $i$ and $n$.  Then $f$ is doubly
symmetric.
\end{thm}

\begin{proof}
Write $f=\sum n_\lambda p_\lambda$.  We assume, by 
way of contradiction, that there is a non-odd $\mu$
with $n_\mu\ne0$, and let $2a$ be the smallest even part of $\mu$.  Set
$$\mathcal{A}=\{\text{non-odd } \lambda \mid m_\lambda\ne0, \mbox{ and the smallest even part of } \lambda \text{ is } 2a\}.$$
We now consider $p(t,-t,x_3,\ldots,x_n)$.  In general,
$$p_k(t,-t,x_3,\ldots,x_n)=\begin{cases}
p_k(x_3,\ldots,x_n),&\text{if $k$ is odd, and}\\
2t^k+p_k(x_3,\ldots,x_n),&\text{if $k$ is even.}
\end{cases}$$
Hence, if $\lambda \in \mathcal{A}$, then
\begin{align*}
p_ \lambda(t,-t,x_3,\ldots,x_n) = & \ p_ \lambda(x_3,\ldots,x_n)+
2t^{2a}p_{\lambda-[2a]}(x_3,\ldots,x_n)\\
&\quad +\mbox{ terms of higher degree in } t,
\end{align*}
where $\lambda-[2a]$ denotes the partition obtained by deleting the part $2a$ from $\lambda$.
The coefficient of $t^{2a}$ in $f(t,-t,x_3,\ldots,x_n)$ is
$$\sum_{\lambda\in\mathcal{A}}n_\lambda p_{\lambda-[2a]}(x_3,\ldots,x_n)=0.$$
Since the shapes $\lambda \in \mathcal{A}$ are distinct, so are the shapes $\lambda-[2a]$, implying that the $p_{\lambda-[2a]}$ are linearly independent.  This contradiction
completes the proof.
\end{proof}

\begin{cor}
Let $f=\sum m_\lambda S_\lambda$ be such that for
all~$n$, $f(t,-t,x_3,\ldots,x_n)$ is not a function of~$t$.  Then
$f$ is doubly symmetric.
\end{cor}

\begin{proof}
Since $f$ is symmetric,
$f(x_1,\ldots,x_n)|_{x_i=-x_j=t}$ does not depend on~$t$ for all
$i,j\le n$, and so the previous theorem applies.
\end{proof}

\end{document}